\documentclass[12pt,twoside]{article}

\date{}

\begin{document}

\centerline{}

\centerline {\Large{\bf Innerproduct Hyperspaces }}

\newcommand{\mvec}[1]{\mbox{\bfseries\itshape #1}}

\centerline{}
\centerline{\bf {Sanjay Roy$^{1}$, T.K. Samanta$^{2}$ }}
\centerline{}
\centerline{$^{1}$Department of Mathematics, South Bantra Ramkrishna Institution, West Bengal, India }
\centerline{E-mail: sanjaypuremath@gmail.com}
\centerline{}
\centerline{$^{2}$Department of Mathematics, Uluberia College, West Bengal, India}
\centerline{E-mail: mumpu$_{-}$tapas5@yahoo.co.in}

\newtheorem{Theorem}{\quad Theorem}[section]

\newtheorem{definition}[Theorem]{\quad Definition}

\newtheorem{theorem}[Theorem]{\quad Theorem}

\newtheorem{remark}[Theorem]{\quad Remark}

\newtheorem{corollary}[Theorem]{\quad Corollary}

\newtheorem{note}[Theorem]{\quad Note}

\newtheorem{lemma}[Theorem]{\quad Lemma}

\newtheorem{result}[Theorem]{\quad Result}

\newtheorem{proposition}[Theorem]{\quad Proposition}

\newtheorem{example}[Theorem]{\quad Example}

{\begin{abstract}
\textbf{\emph{The main purpose of this paper is to generalize and develop a few basic properties of the innerproduct space on a hypervector space. On this  hypervector space we define the norm. Also we establish a important relation between normed hyperspaces and innerproduct hyperspaces.}}
\end{abstract}

\section{Introduction}
\smallskip\hspace{.5cm}The algebraic hyperstructure is a generalization of the concept of algebraic structure. In the usual algebraic structure the composition of two elements is a single element. But in the  algebraic hyperstructure the composition of two elements is a set. The hyperstructure was first introduced by F. Marty \cite{Marty1} in 1934. Then he established the definition of hypergroup \cite{Marty2} in 1935. M. Krasner \cite{Krasner} , a great researcher in this area, introduced the notions of hyperring and hyperfield . Since then many researchers ( \cite{Corsini}, \cite{Marty3} ) have studied and developed the concept of different types of hyperstructures in different views.  Then M. S. Tallini ( \cite{Tallini1}, \cite{Tallini2}, \cite{Tallini3} ) introduced the notion of  hypervector spaces. Then in 2005 R. Ameri \cite{Ameri} also studied this spaces extensively. In our previous paper \cite{Roy} we also introduced the notion of a hypervector spaces in more generalized form than the previous concept of hypervector space and thereafter established a few useful theorems in this space.\\
\smallskip\hspace{.5cm}In this paper we recall first the main  definitions and results of the hyperstructure which will be needed. Next in $\S3$, we present the more generalized definition of the hypervector space than the previous definitions and also we present a few definitions and state a usual theorem on this hypervector space.\\
\smallskip\hspace{.5cm}In $\S4$ we first introduce the definition of normed hyperspaces and then we establish the definition of innerproduct hyperspaces. Also here we establish the relation between the structures of norm and innerproduct on hyperspaces. Next we present the concept of orthogonal Set and the orthonormal set. lastly it is seen that every finite dimensional innerproduct hyperspace possesses an orthogonal basis.

\section{Preliminaries}
\smallskip\hspace{.5cm}This section contain some basic definitions and preliminary results which will be needed.

\begin{definition} \cite{Nakassis}
A hyperoperation over a non-empty set X is a mapping of $X\times X$ into the set of all non-empty subsets of X.
\end{definition}
\begin{definition}\cite{Nakassis}
A non-empty set X with exactly one hyperoperation $`\#$' is called a hypergroupoid.\\
Let $(X\, ,\, \#)$ be a hypergroupoid. For every point $x\in X$ and every non-empty subset A of X , we defined  $x\;\#\; A\;=\;\bigcup_{a\,\in\,A}\{ x \;\#\; a\}$.
 \end{definition}
\begin{definition}\cite{Nakassis}
A hypergroupoid $(X \,,\,\#)$ is called a semihypergroup if $x\#(y\#z)=(x\#y)\#z$  for all $x,\,y,\,z \,\in\, X$.
\end{definition}
\begin{definition}\cite{Nakassis}
A hypergroupoid $(X \,,\,\#)$ is called a hypergroup if\\
$( i )$\hspace{0.5cm}   $a\#(b\#c)=(a\#b)\#c$ for all $a,\,b,\,c \,\in\, X$.\\
$( ii )$\hspace{0.4cm}  $\exists\, 0\in\, X$ such that for every a $\in\, X$, there is unique element $b \in X$ for which $0 \,\in a\# b$ and $0 \,\in b \# a$. Here $b$ is denoted by $-a$.\\
$( iii )$ \hspace{0.3cm} For all $a , b , c \in X$ if $a\in b \# c$, then $b\in a \#(- c)$.
\end{definition}
\begin{proposition}\cite{Nakassis}
$(i)$ In a hypergroup $(X \,,\, \#)$, $-(- a) = a$, $\forall  a \in X$.

$(ii)$ $0\;\#\; a=\{a\}$, $\forall a \in X$ if $(X\,,\,\#)$ is a commutative hypergroup.

$(iii)$ In a commutative hypergroup $($X , $\#)$, $0$ is unique.
\end{proposition}
\begin{note}\cite{Nakassis}
 In a hypergroup, if the element $0$ is unique, then $0$ is called the zero element of the hypergroup and $b$ is called the additive inverse of $a$ if $0\in a\#b$ and $0\in b\# a$.
\end{note}
\begin{definition}\cite{Roy}
A hyperring is a non-empty set equipped with a hyperaddition $`\#$' and a multiplication $`.$' such that $(X\,,\,\#)$ is a commutative hypergroup and $(X\,,.)$ is a semigroup and the multiplication is distributive across the hyperaddition both from the left and from the right and $a.0 = 0.a = 0$,  $\forall\, a\,\in X$, where $0$ is the zero element of the hyperring.
\end{definition}
\begin{definition}\cite{Roy}
A hyperfield is a non-empty set X equipped with a hyperaddition $`\#$' and a multiplication $`.$' such that\\
$(\,i\,)$\hspace{0.5cm}$(X\,,\;\#\;,\,.)$ is a hyperring.\\
$(\,ii\,)$\hspace{0.4cm}$\exists$ an element $1 \in$ X, called the identity element such that $a.1 = a$, $\forall \,a\,\in X.$ \\
$(\,iii\,)$\hspace{0.3cm}For each non-zero element $a\in X$, $\exists$ an element $a^{-1}\in X$ such that $a.a^{-1}=1.$\\
$(\,iv\,)$\hspace{0.4cm} $a.b = b.a,$ $\forall\, a, b\in X$.
\end{definition}

\section{Hypervector spaces}
\smallskip\hspace{.5cm}In this section we modify some definitions in hypervector space.
\begin{definition}
Let $(F\,,\oplus\,,.)$ be a hyperfield and $(V\,,\,\#)$ be an additive commutative hypergroup. Then $V$ is said to be a hypervector space or hyperspace over the hyperfield $F$ if there exists a hyperoperation\\
 $\ast: \;F \;\times\, V \,\rightarrow\;P^{\ast}\,(\,V\,)$ such that\\
$(\,i\,)$\hspace{0.5cm}a $\ast\,(\alpha\;\#\;\beta)\subseteq a\ast\alpha\;\,\#\;\,a\ast\beta,$  \hspace{1cm}$\forall\, a\in F\;and\;\forall\, \alpha\, , \beta \in$ $V$.\\
$(\,ii\,)$\hspace{0.4cm}$(a\,\oplus\,b)\ast\alpha \subseteq a\ast\alpha\;\#\;b\ast\alpha,$   \hspace{1cm} $\forall\,a , b\in F\;and\;\forall\,\alpha\in
$ $V$.\\
$(\,iii\,)$\hspace{0.3cm}$(a\,.\,b)\ast\alpha=a\ast(b\ast\alpha),$   \hspace{2cm} $\forall\, a , b \in$F and $\forall\, \alpha\in$ $V$.\\
$(\,iv\,)$\hspace{0.4cm}$(-a)\ast\alpha=a\ast(-\alpha), \hspace{2.7cm}$  $\forall\, a\in$ $F$ and $\forall\, \alpha\in$ V.\\
$(\,v\,)$\hspace{0.5cm}$\alpha\in 1_{F}\ast\alpha$, $\theta\in 0 \ast\alpha$ and $0\ast\theta\,=\,\theta,$\hspace{1cm}$\forall\, \alpha\in$ $V$.\\
 where $1_{F}$ is the identity element of $F,$  $0$ is the zero element of $F$ and $\theta$ is zero vector of $V$ and $P^{\ast}(\,V\,)$ is the set of all non-empty subset of $V$.
 \end{definition}

\begin{result}
 In a hypervector space $V,$ Show that $-\alpha\in -1_{F}\ast\alpha.$
\end{result}
 \textbf{Proof:} Since  $\alpha\in 1_{F}\ast\alpha$, $\forall$ $\alpha\in V$
 $\Rightarrow-\alpha\in 1_{F}\ast-\alpha$, $\forall$ $\alpha\in V$\\
 $\Rightarrow-\alpha\in -1_{F}\ast\alpha$, $\forall$ $\alpha\in V,$ as $(-1_{F})\ast\alpha\,=\,1_{F}\ast(-\alpha).$

\begin{definition}
Let $(V,\,\#,\,\ast,\,F)$ be a hypervector space. A subset $A=\{\alpha_{\lambda}\}_{\lambda\in\Lambda}$ of $V$ is set to a linearly dependent set if there exists a finite subset $\{\alpha_{1},\alpha_{2},.\,.\,\,.,\alpha_{n}\}\subseteq\,$ A such that
$\theta\,\in\,\lambda_{1}\ast \alpha_{1}\,\#\,\lambda_{2}\ast
\alpha_{2}\,\#\,.\,.\,.\,\#\,\lambda_{n}\ast \alpha_{n}$ for some $\lambda_{1},\,\lambda_{2},\,.\,.\,.\,,\,
\lambda_{n}\,$(not all zeros)$\in$ F.\\
\smallskip\hspace{.5cm}Otherwise A is said to a linearly independent set.
\end{definition}
\begin{definition}
A subset A of a hypervector space $(V,\,\#,\,\ast,\,F)$ is said to be weak linearly independent set if for all finite subset $\{\alpha_{1},\alpha_{2},.\,.\,\,.,\alpha_{n}\}$ of A and $0\ast P\,=\,\lambda_{1}\ast \alpha_{1}\,\#\,\lambda_{2}\ast \alpha_{2}\,\#\,.\,.\,.\,\#\,\lambda_{n}\ast \alpha_{n},$ for some $\lambda_{1},\,\lambda_{2},\,.\,.\,.\,,\,
\lambda_{n}\,\in$ F and for some $P\,\subseteq\,V \Rightarrow\, \lambda_{1}=\lambda_{2}=...=\lambda_{n}=0.$\\
\smallskip\hspace{.5cm}Otherwise A is said to a weak linearly dependent set.
\end{definition}
\begin{definition}
A linearly independent subset A of a hypervector space\\
$(V,\,\#,\,\ast,\,F)$ is called a basis if for every $\alpha\in V,$ there are n $(\in \mathbf{N})$ elements $\lambda_{1},\,\lambda_{2},\,.\,.\,.\,,\lambda_{n}\in F$ and $\alpha_{1},\alpha_{2},\,.\,.\,.\,,\alpha_{n}\in$ A such that\\
$ \alpha\,\in\,\lambda_{1}\ast \alpha_{1}\,\#\,\lambda_{2}\ast \alpha_{2}\,\#\,.\,.\,.\,\#\,\lambda_{n}\ast \alpha_{n}.$\\
The hypervector space is said to be finite dimensional if it has a finite basis.
\end{definition}
\begin{theorem}
Let $(V,\,\#,\,\ast,\,F)$ be a finite dimensional hypervector space and  $\{\alpha_{1},\alpha_{2},\,.\,.\,\,.,\,\alpha_{n}\}\subseteq\,$ V be a basis of V. If
$ \alpha\,\in\,\lambda_{1}\ast \alpha_{1}\,\#\,\lambda_{2}\ast \alpha_{2}\,\#\,.\,.\,.\,\#\,\lambda_{n}\ast \alpha_{n}$ for some $\lambda_{1},\,\lambda_{2},\,.\,.\,.,\,\lambda_{n}\in F$ and for some element $\alpha\in V$. Then the set $\{\alpha,\,\alpha_{1},\alpha_{2},\,.\,.\,\,.,\,\alpha_{n}\}$ is linearly dependent.
\end{theorem}
\textbf{proof:} obvious.

\section{Innerproduct on hyperspaces}
\smallskip\hspace{.5cm}In this section at first we define a norm on a hypervector space. Then we introduce the definition of innerproduct hyperspace and deduce some important theorems.
\begin{definition}
Let $\mathbf{R}$ be the set of all real numbers. The hyperfield defined on $\mathbf{R}$ is called the real hyperfield.
\end{definition}

\begin{definition}
Let $(V\,,\#\,,\ast)$ be a hypervector space over the real hyperfield $\mathbf{R}$. A norm on V is a mapping  $\|\cdot\|:V\rightarrow \mathrm{R},$ where $\mathrm{R}$ is a usual real space, such that for all $a\in\mathbf{R}$ and $\alpha, \beta\in$ V conditions hold\\
$(i)$ $\|\alpha\|\geq 0$.\\
$(ii)$ $\|\alpha\|=0$ if and only if $\alpha=\theta.$\\
$(iii) \sup\|\alpha\#\beta\|\leq \|\alpha\| + \|\beta\|,\;\;\;$where $\|\alpha\,\#\,\beta\|\,=\,\{\,\|x\|,\;\;x\,\in\,\alpha\,\#\,\beta \,\}.$\\
$(iv) \sup\|a\ast \alpha\|\,\leq \,|a|\|\alpha\|,\;\;$where $\|a\ast\alpha\|\,=\,\{\,\|x\|,\;\;x\,\in\,a\ast \alpha \,\}.$
\end{definition}

If $\|\cdot\|$ is a norm on $V$ then $(V, \#, \ast)$ is said to be a normed hypervector space or normed hyperspace.\\

\begin{definition}
Let $(V,\,\#,\,\ast)$ be a hypervector space over the real hyperfield $\mathbf{R}$. By an innerproduct on V, we mean a mapping
$<\cdot\,,\,\cdot>\, : V\times V\rightarrow \mathrm{R},$ where $\mathrm{R}$ is a usual real space, such that\\
$(i)$ $<\alpha\,,\,\alpha>\;\;> 0\;,\forall\,\alpha\,(\neq\,\theta)\in$ V.\\
$(ii)$ $<\alpha\,,\,\alpha>\;=\,0$ iff $\alpha=\theta.$\\
$(iii)$ $<\alpha\,,\,\beta>\;=\;<\beta\,,\,\alpha>\;\;\forall\, \alpha,\,\beta\,\in$ V.\\
$(iv)$ $ \sup<\alpha\#\beta\,,\,\gamma>\;=\,<\alpha\,,\gamma>\,+\,<\beta\,,\gamma>\;\;\forall\, \alpha,\, \beta,\, \gamma\,\in\, V,$\\
\smallskip\hspace{4cm}Where $<\alpha\#\beta\,,\,\gamma>\;=\,\{<\delta\,,\,\gamma>\;:\;\delta \in \alpha\#\beta\}.$\\
$(v)$ $\sup<a\ast\alpha\,,\,\beta>\;=\,a<\alpha\,,\,\beta>\;\;\forall\, \alpha,\, \beta\,\in\, $V and $\forall\, a\,\in\, \mathbf{R}$.
\end{definition}

\begin{definition}
A hypervector space together with an innerproduct is called an innerproduct hyperspace.
\end{definition}

\begin{result}
Let $(V,\,\#,\,\ast)$ be a hypervector space over the real hyperfield $\mathbf{R}$ and $\alpha, \beta, \gamma\,\in$ V. Then the following statements are true\\
$(i)$  $\sup<\alpha\,,\,\beta\#\gamma>\,=\,<\alpha\,,\beta>\,+\,<\alpha\,,\gamma>.$\\
$(ii)$ $\sup<\alpha\,,\,a\ast\beta>\,=\,a<\alpha\,,\,\beta>$.\\
$(iii)$ $<\alpha\,,\,\theta>\,=\,<\theta\,,\,\alpha>\,=\,0.$
\end{result}
\textbf{proof:}$\;$ Obvious.

\begin{lemma}\label{l1}
Let V be an innerproduct hyperspace over the real hyperfield $\mathbf{R}$, then
$\;\sup<\alpha\,\#\,a\ast\beta\,,\,\gamma\,\#\,b\ast\delta>\,=\,<\alpha\,,\,\gamma>+a <\beta\,,\,\gamma>+b\,<\alpha\,,\,\delta>\\
\smallskip\hspace{9cm}+a b <\beta\,,\,\delta>$\\
where $\alpha, \beta, \gamma, \delta \in$ V; a, b $\in \mathbf{R}$ and
\\$<\alpha\,\#\,a\ast\beta\,,\,\gamma\,\#\,b\ast\delta>\,=\,\{\,<x\,,\,y>\,,\;\;x\,\in \alpha\,\#\,a\ast\beta,\,y\,\in \gamma\,\#\,b\ast\delta\,\}.$
\end{lemma}
\textbf{proof:}
$\;\sup<\alpha\,\#\,a\ast\beta\,,\,\gamma\,\#\,b\ast\delta>$\\
= $\sup\{\sup<\alpha\,\#\,x\,,\,y>,\;x\in a\ast\beta,\;\; y\in\gamma\,\#\,b\ast\delta\}$\\
= $\sup\{<\alpha\,,\,y>\, + \,<x\,,\,y>,\;x\in\, a\ast\beta,\;\; y\,\in\,\gamma\,\#\,b\ast\delta\}$\\
= $\sup\{<\alpha\,,\,y>\, + \,\sup<a\ast\beta\,,\,y>,\; \;y\,\in\,\gamma\,\#\,b\ast\delta\}$\\
= $\sup\{<\alpha\,,\,y>\, + \,a <\beta\,,\,y>,\;\; y\,\in\,\gamma\,\#\,b\ast\delta\}$\\
= $\sup\{\sup<\alpha\,,\,\gamma\,\#\,\eta> \,+\, a \sup<\beta\,,\,\gamma\,\#\,\eta>,\;\;\eta\,\in\, b\ast\delta\}$\\
= $\sup\{<\alpha\,,\,\gamma> + <\alpha\,,\,\eta> \,+\, a <\beta\,,\,\gamma>\, + \,a <\beta\,,\,\eta>,\;\;\eta\,\in \,b\ast\delta\}$\\
= $<\alpha\,,\,\gamma>\, +\, \sup<\alpha\,,\,b\ast\delta>\, + \,a <\beta\,,\,\gamma>\, +\, a \sup<\beta\,,\,b\ast\delta>$\\
= $<\alpha\,,\,\gamma> \,+\, b <\alpha\,,\,\delta> \,+\, a <\beta\,,\,\gamma>\, +\, ab <\beta\,,\,\delta>.$

\begin{result}\label{r1}
Let V be a hypervector space over the real hyperfield $\mathbf{R}$ and \\
$<\cdot , \cdot>$ be an innerproduct on V. We define a mapping f : $V\rightarrow \mathrm{R}$  $(\mathrm{R}$ is a usual real space$)$ by \\
$f(\alpha)= \sqrt{<\alpha\,,\,\alpha>}$, then $\mid<\alpha\,,\,\beta>\mid \,\leq f(\alpha) f(\beta)\;\;\forall\, \alpha, \beta \in V$.
\end{result}
\textbf{proof:}$\;$ If $\beta = \theta$, then $<\alpha\,,\,\theta> = 0$ i.e $\mid<\alpha\,,\,\beta>\mid = 0$ for $\beta = \theta.$\\
Again $f(\beta) = \sqrt{<\beta\,,\,\beta>} = 0$ for $\beta = \theta.$\\
Therefore $f(\alpha)f(\beta) = 0.$\\
Hence $\mid<\alpha\,,\,\beta>\mid = f(\alpha) f(\beta)$ if $\beta=\theta.$\\
Similarly, this holds for $\alpha = \theta$.\\
So we assume that $\alpha, \beta \neq \theta$, then for every scalar $a \in \mathbf{R}$ we have\\
$\sup<\alpha \,\# \,a\ast\beta\,,\,\alpha\, \#\, a\ast\beta>\;\;\geq 0$\\
$\Rightarrow\,<\alpha\,,\,\alpha>\, +\, a<\beta\,,\,\alpha>\, + \,a<\alpha\,,\,\beta>\, +\, a^{2}<\beta\,,\,\beta>\;\;\geq 0\;\;$\\
\smallskip\hspace{8cm}[ by lemma \ref{l1}]\\
$\Rightarrow\,<\alpha\,,\,\alpha> \,+ \,a<\alpha\,,\,\beta>\, +\, a<\alpha\,,\,\beta>\, +\, a^{2}<\beta\,,\,\beta>\;\;\geq 0$\\
$\Rightarrow\,<\alpha\,,\,\alpha>\, +\, 2a<\alpha\,,\,\beta>\, + \,a^{2}<\beta\,,\,\beta>\;\;\geq 0$\\
$\Rightarrow\, (\,f(\alpha)\,)^{2} \,+\, 2a<\alpha\,,\,\beta>\, + \,a^{2}(\,f(\beta)\,)^{2}\;\;\geq 0$\\
put $a \,=\, -\frac{<\alpha\,,\,\beta>}{(\,f(\beta)\,)^{2}}$\\
Therefore $(\,f(\alpha)\,)^{2}\, -\, 2\frac{<\alpha\,,\,\beta>}{(\,f(\beta)\,)^{2}}<\alpha\,,\,\beta> \,+\, \{-\frac{<\alpha\,,\,\beta>}{(\,f(\beta)\,)^{2}}\}^{2}(\,f(\beta)\,)^{2}\;\;\geq 0$\\
$\Rightarrow\,(\,f(\alpha)\,)^{2}\, - \,2\frac{(\,<\alpha\,,\,\beta>\,)^{2}}{(\,f(\beta)\,)^{2}} \,+ \, \frac{(\,<\alpha\,,\,\beta>\,)^{2}}{(\,f(\beta)\,)^{2}}\;\;\geq 0$\\
$\Rightarrow\, (\,<\alpha\,,\,\beta>\,)^{2}\,\leq\,(\,f(\alpha)\,)^{2} (\,f(\beta)\,)^{2}$\\
$\Rightarrow\, \mid<\alpha\,,\,\beta>\mid\,\leq\,f(\alpha) f(\beta).$

\begin{result}
The function f defined in the result $\ref{r1}$ satisfy the following properties\\
$(i)$ $f(\alpha)\;\geq 0.$\\
$(ii)$ $f(\alpha)\, =\, 0$  iff $\alpha = \theta$.\\
$(iii)$ $\sup f(a\ast\alpha) \leq |a| f(\alpha),\;\;\forall\, \alpha\, \in$ V and $\forall\, a \,\in \mathbf{R},\;\;$\\where $f(a\ast\alpha)\,=\,\{\,f(x),\;\;x\,\in\,a\ast\alpha\,\}.$\\
$(iv)$ $\sup f(\alpha \# \beta)\,\leq\, f(\alpha)\, + \,f(\beta),\;\;$where $f(\alpha \# \beta)\,=\,\{\,f(x),\;\;x\,\in\,\alpha\,\#\,\beta\,\}.$
\end{result}
\textbf{proof:}$\;$ Prove of $(i)$ and $(ii)$ are obvious.\\
(iii) $\sup f(\,a\ast\alpha\,)\,=\,\sup\{\,f(x)\;:\;\;x\in\,a\ast\alpha\,\}$\\
 \smallskip\hspace{3.3cm}= $\sup\{\,\sqrt{<x\,,\,x>}\;:\;\;x\in\,a\ast\alpha\,\}$\\
\smallskip\hspace{3.3cm}= $\sqrt{\sup\{<x\,,\,x>\;:\;\;x\in\, a\ast\alpha\}}$\\
\smallskip\hspace{3.3cm} $\leq \sqrt{\sup<a\ast\alpha\,,\,a\ast\alpha>}$\\
\smallskip\hspace{3.3cm}= $\sqrt{a^{2} <\alpha\,,\,\alpha>}$\\
\smallskip\hspace{3.3cm}= $|a|\sqrt{<\alpha\,,\,\alpha>}$\\
\smallskip\hspace{3.3cm}= $|a| f(\alpha).$\\
Therefore  $\sup f(\,a\ast\alpha\,)\leq |a| f(\alpha).$\\

\textbf{(iv)} $\sup f(\alpha\,\#\,\beta)\,=\, \sup\{\,f(x)\;:\;\;x\in\,\alpha\,\#\,\beta\,\}$\\
\smallskip\hspace{3.9cm}= $\sup\{\,\sqrt{<x\,,\,x>}\;:\;\;x\in\,\alpha\,\#\,\beta\,\}$\\
\smallskip\hspace{3.9cm}= $\sqrt{\sup\{\,<x\,,\,x>\;:\;\;x\in\,\alpha\,\#\,\beta\,\}}$\\
\smallskip\hspace{3.9cm} $\leq \sqrt{\sup<\alpha\,\#\,\beta\,,\,\alpha\,\#\,\beta>}$\\
\smallskip\hspace{3.9cm}= $\sqrt{<\alpha\,,\,\alpha>\, +\, <\beta\,,\,\beta>\, +\, <\alpha\,,\,\beta>\, +\, <\beta\,,\,\alpha>}$\\
\smallskip\hspace{3.9cm}= $\sqrt{<\alpha\,,\,\alpha>\, + \,<\beta\,,\,\beta> \,+ \,2<\alpha\,,\,\beta>}$\\
\smallskip\hspace{3.9cm} $\leq \,\sqrt{(\,f(\alpha)\,)^{2}\, + \,(\,f(\beta)\,)^{2}\, +\, 2 f(\alpha) f(\beta)},\;$ $[$by result $\ref{r1}]$\\
\smallskip\hspace{3.9cm}= $f(\alpha)\, + \,f(\beta).$\\
Therefore $\;\; \sup f(\alpha\,\#\,\beta)\leq f(\alpha)\, + \,f(\beta).$

\begin{remark}
Observing the properties $(i),\,(ii),\,(iii)$ and $(iv)$ we can say that the function f satisfies all the conditions of norm and hence f is a norm on the hypervector space V. So we can replace f by $\parallel\cdot\parallel$.\\
So every innerproduct hyperspace is a normed hyperspace.
 \end{remark}

\begin{result}
If $\alpha,\,\beta$ be any two vectors in an innerproduct hyperspace V, then
 \hspace{.3cm}$\sup\|\alpha\,\#\,\beta\|^{2}\,+\,\sup\|\alpha\,\#\,(-\beta)\|^{2}\,\leq\,2\|\alpha\|^{2}\,+\,2\|\beta\|^{2}.$
\end{result}
\textbf{proof:}$\;\sup\|\alpha\,\#\,\beta\|^{2}\,=\,\sup\{\,\|x\|^{2}\;:\;\;x\in\,\alpha\,\#\,\beta\,\}$\\
\smallskip\hspace{3.8cm}= $\,\sup\{\,<x\,,\,x>\;:\;\;x\in\,\alpha\,\#\,\beta\,\}$\\
\smallskip\hspace{3.8cm} $\leq\,\sup<\alpha\,\#\,\beta\,,\,\alpha\,\#\,\beta>$\\
\smallskip\hspace{3.8cm}= $\,<\alpha\,,\,\alpha>\,+\,<\beta\,,\,\alpha>\,+\,<\alpha\,,\,\beta>\,+\,<\beta\,,\,\beta>$\\
\smallskip\hspace{3.8cm}= $\,\|\alpha\|^{2}\,+\,2<\alpha\,,\,\beta>\,+\,\|\beta\|^{2}.$\\
Therefore $\;\sup\|\alpha\,\#\,\beta\|^{2}\,\leq\,\|\alpha\|^{2}\,+\,2<\alpha\,,\,\beta>\,+\,\|\beta\|^{2}.$\\

$\sup\|\alpha\,\#\,(-\beta)\|^{2}\\
=\,\sup\{\,\|x\|^{2}\;:\;\;x\,\in\,\alpha\,\#\,(-\beta)\,\}$\\
= $\,\sup\{\,<x\,,\,x>\;:\;\;x\,\in\,\alpha\,\#\,(-\beta)\,\}$\\
$\leq\,\sup<\alpha\,\#\,(-\beta)\,,\,\alpha\,\#\,(-\beta)>$\\
 =$\,<\alpha\,,\,\alpha>\,+\,<-\beta\,,\,\alpha>\,+\,<\alpha\,,\,-\beta>\,+\,<-\beta\,,\,-\beta>$\\
$\leq\; \,<\alpha\,,\,\alpha>\,+\,\sup<-1\ast\beta\,,\,\alpha>\,+\,\sup<\alpha\,,\,-1\ast\beta>\,+\\$
\smallskip\hspace{1cm}$ \sup<-1\ast\beta\,,\,-1\ast\beta>$, as $-\alpha\in-1\ast\alpha, \forall\,\alpha\in V$\\
= $\,\|\alpha\|^{2}\,-\,<\beta\,,\,\alpha>\,-\,<\alpha\,,\,\beta>\,+\,\|\beta\|^{2}$\\
= $\,\|\alpha\|^{2}\,-\,2<\alpha\,,\,\beta>\,+\,\|\beta\|^{2}.$\\
Therefore $\sup\|\alpha\,\#\,(-\beta)\|^{2}\,\leq\,\|\alpha\|^{2}\,-\,2<\alpha\,,\,\beta>\,+\,\|\beta\|^{2}.$\\
Hence $\sup\|\alpha\,\#\,\beta\|^{2}\,+\,\sup\|\alpha\,\#\,(-\beta)\|^{2}\,\leq\,2\|\alpha\|^{2}\,+\,2\|\beta\|^{2}.$

\begin{definition}
Let V be an innerproduct hyperspace. The vectors $\alpha$ is said to be orthogonal to $\beta$ in V if $<\alpha\,,\,\beta> = 0$. A subset S of V is orthogonal if any two distinct vectors in S are orthogonal.\\
\smallskip\hspace{.5cm}A vector $\alpha$ in V is said to be unit vector if $\|\alpha \| = 1.$\\
\smallskip\hspace{.5cm}A subset S of V is orthonormal if S is orthogonal and consists entirely of unit vectors.
\end{definition}

\begin{definition}
An orthogonal $($orthonormal$)$ subset A of V is called a orthogonal basis $($orthonormal basis$)$ if for every $\alpha\in$ V, there are n $(\in \mathbf{N})$ elements $\lambda_{1},\,\lambda_{2},\,.\,.\,.\,,\lambda_{n}\in$ F and $\alpha_{1},\alpha_{2},\,.\,.\,.\,,\alpha_{n}\in$ A such that\\
$ \alpha\,\in\,\lambda_{1}\ast \alpha_{1}\,\#\,\lambda_{2}\ast \alpha_{2}\,\#\,.\,.\,.\,\#\,\lambda_{n}\ast \alpha_{n}.$
\end{definition}

\begin{theorem}
Let $(\,V\,,\,<\cdot\,,\,\cdot>\,)$ be an innerproduct hyperspace. Then any non-null orthogonal subset of V is weak linearly independent.
\end{theorem}
\textbf{proof:}$\;$
Let S be any non-null orthogonal subset of V and $\alpha_{1},\,\alpha_{2},\,.\,.\,.,\,\alpha_{n}$ be any finite vectors of S.
Let $0\ast P = \lambda_{1}\ast\alpha_{1}\,\#\,\lambda_{2}\ast\alpha_{2}\,\#\,.\,.\,.\,\lambda_{n}\ast\alpha_{n}\;\;$ for some $P \subseteq V$ and $\lambda_{1},\,\lambda_{2},\,.\,.\,.,\,\lambda_{n} \in \mathbf{R}$.\\
$\Rightarrow\, <0\ast P\,,\alpha_{i}> = <\lambda_{1}\ast\alpha_{1}\,\#\,\lambda_{2}\ast\alpha_{2}\,\#\,.\,.\,.\,\#\,\lambda_{n}\ast\alpha_{n}\,,\alpha_{i}>$, $\;\;\;i\,\in\,\{1, 2, . . .,n\}$\\
$\Rightarrow\, \sup<0\ast P\,,\alpha_{i}>\, = \,\sup<\lambda_{1}\ast\alpha_{1}\,\#\,\lambda_{2}\ast\alpha_{2}\,\#\,.\,.\,.\,\lambda_{n}\ast\alpha_{n}\,,\alpha_{i}>,$ \\
\smallskip\hspace{8cm}$i\in\{1, 2, . . .,n\}$\\
$\Rightarrow\, 0\sup<P\,,\,\alpha_{i}>\, = \,\lambda_{1}<\alpha_{1}\,,\,\alpha_{i}> + \lambda_{2}<\alpha_{2}\,,\,\alpha_{i}> + \,.\,.\,.\,+ \lambda_{n}<\alpha_{n}\,,\,\alpha_{i}>$\\
$\Rightarrow\, 0\, = \,\lambda_{i}<\alpha_{i}\,,\,\alpha_{i}>$\\
$\Rightarrow \,\lambda_{i}\, =\, 0.$\\
Therefore $\lambda_{1} \,=\, \lambda_{2}=\,.\,.\,.\,=\lambda_{n} = 0$.\\
Therefore $\{\alpha_{1},\,\alpha_{2},\,.\,.\,.\,\alpha_{n}\}$ is weak linearly independent.\\
Hence any non-null orthogonal subset of V is weak linearly independent.

\begin{theorem}
Let $(V\,,\#\,,\ast)$ be an innerproduct hyperspace over the real hyperfield $\mathbf{R}$ and $S= \{\alpha_{1},\,\alpha_{2},\,.\,.\,.,\,\alpha_{n}\}$ be an orthogonal subset of V consisting of non-zero vectors. If a vector $\alpha\,\in\,$V can be expressed as\\
$1\ast \alpha\,=\,\lambda_{1}\ast \alpha_{1}\,\#\,\lambda_{2}\ast \alpha_{2}\,\#\,.\,.\,.\,\#\,\lambda_{n}\ast \alpha_{n},\;\;$
for some $\lambda_{1},\,\lambda_{2},\,.\,.\,.\,,\lambda_{n}\in \mathbf{R}$ then\\
\[1\ast\alpha\,=\,\mathop {\sum }\limits_{i=1 }\limits^{n}\frac{<\alpha\,,\,\alpha_{i}>}{\|\alpha_{i}\|^{2}}\ast\alpha_{i}\]
\end{theorem}
\textbf{proof:}$\;$ Let $1\ast \alpha=\lambda_{1}\ast \alpha_{1}\,\#\,\lambda_{2}\ast \alpha_{2}\,\#\,.\,.\,.\,\#\,\lambda_{n}\ast \alpha_{n},$
for some $\lambda_{1},\,\lambda_{2},\,.\,.\,.\,,\lambda_{n}\in \mathbf{R}$.\\
Then for $1\leq j \leq n$, we have\\
$<1\ast \alpha\,,\,\alpha_{j}>\,=\,<\lambda_{1}\ast \alpha_{1}\,\#\,\lambda_{2}\ast \alpha_{2}\,\#\,.\,.\,.\,\#\,\lambda_{n}\ast \alpha_{n}\,,\,\;\;\alpha_{j}>\;\;\;$\\
$\Rightarrow\,\sup<1\ast \alpha\,,\,\alpha_{j}>\,=\,\sup<\lambda_{1}\ast \alpha_{1}\,\#\,\lambda_{2}\ast \alpha_{2}\,\#\,.\,.\,.\,\#\,\lambda_{n}\ast \alpha_{n}\,,\,\;\;\alpha_{j}>$\\
$\Rightarrow\, 1<\alpha\,,\,\alpha_{j}>\,=\,\lambda_{1}<\alpha_{1}\,,\,\alpha_{j}>\,+\,\lambda_{2}<\alpha_{2}\,,\,\alpha_{j}>\,+\,.\,.\,.\,+\,\lambda_{n}
<\alpha_{n}\,,\,\alpha_{j}>$\\
$\Rightarrow\,<\alpha\,,\,\alpha_{j}>\,=\,\lambda_{j}<\alpha_{j}\,,\,\alpha_{j}>$\\
\smallskip\hspace{2.4cm}= $\lambda_{j}\|\alpha_{j}\|^{2}$\\
$\Rightarrow\, \lambda_{j}\,=\,\frac{<\alpha\,,\,\alpha_{j}>}{\|\alpha_{j}\|^{2}},\;\;\;j=1, 2, . . ., n$\\
$\Rightarrow\,1\ast \alpha\,=\,\frac{<\alpha\,,\,\alpha_{1}>}{\|\alpha_{1}\|^{2}}\ast \alpha_{1}\,\#\,\frac{<\alpha\,,\,\alpha_{2}>}{\|\alpha_{2}\|^{2}}\ast \alpha_{2}\,\#\,.\,.\,.\,\#\,\frac{<\alpha\,,\,\alpha_{n}>}{\|\alpha_{n}\|^{2}}\ast \alpha_{n}$\\
Therefore \[1\ast\alpha\,=\,\mathop {\sum }\limits_{i=1 }\limits^{n}\frac{<\alpha\,,\,\alpha_{i}>}{\|\alpha_{i}\|^{2}}\ast\alpha_{i}.\]

\begin{theorem}\label{t1}
Let $(V\,,\#\,,\ast)$ be an innerproduct hyperspace over the real hyperfield $\mathbf{R}$ and $S=\{w_{1},\,w_{2},\,.\,.\,.\,,w_{n}\}$ be a linearly independent subset of V. Define $S^{'}=\{v_{1},\,v_{2},\,.\,.\,.,\,v_{n}\}$, where $v_{1}=w_{1}$ and \\
\[v_{k}\,\in\, w_{k}\,\#\,- \mathop {\sum }\limits_{j=1 }\limits^{k-1}\frac{<w_{k}\,,\,v_{j}>}{\|v_{j}\|^{2}}\ast v_{j},\, for\, 2 \leq k \leq n\;\;\;\cdots\cdots\cdots (1)\]\\
Then $S^{'}$ is an orthogonal set of non-zero vectors such that $span(S^{'})\,=\,span(S)$
\end{theorem}
\textbf{proof:}$\;$
We prove the theorem by the mathematical induction on n, the number of vectors in S.\\
For $k = 1, 2, . . ., n,$ let $S_{K}\,=\,\{w_{1},\,w_{2},\,.\,.\,.,\,w_{k}\}$ and $S_{K}^{'}\,=\,\{v_{1},\,v_{2},\,.\,.\,.,\,v_{k}\}$.\\
If n=1, then the theorem is proved by taking $S^{'}_{1}\,=\,S_{1}$.\\
Assume then that the set $S^{'}_{K-1}\,=\,\{v_{1},\,v_{2},\,.\,.\,.,\,v_{k-1}\}$ with the desired properties has been constructed by the repeated use of $(1)$. We show that the  set $S^{'}_{K}\,=\,\{v_{1},\,v_{2},\,.\,.\,.,\,v_{k-1},\,v_{k}\}$ also has the desired properties, where $v_{k}$ is obtained from $S^{'}_{k-1}$ by $(1)$.\\
We choose $v_{k}$ from $(1)$ such that $v_{k}\neq \theta$\\
If \[\{\theta\}\,=\, w_{k}\,\#\,-\mathop {\sum }\limits_{j=1 }\limits^{k-1}\frac{<w_{k}\,,\,v_{j}>}{\|v_{j}\|^{2}}\ast v_{j}\]
Then \[ w_{k}\,\in\,\mathop {\sum }\limits_{j=1 }\limits^{k-1}\frac{<w_{k}\,,\,v_{j}>}{\|v_{j}\|^{2}}\ast v_{j}\]
$\Rightarrow w_{k}\,\in\,span(S^{'}_{k-1})\,=span(S_{k-1})$, which contradicts the assumption that $S_{K}$ is linearly independent.\\
For $1\,\leq\,i\,\leq\,k-1$, it follows from $(1)$ that\\
\[<v_{k}\,,\,v_{i}>\,\;\;\in\;\;\,< w_{k}\,\#\,- \mathop {\sum }\limits_{j=1 }\limits^{k-1}\frac{<w_{k}\,,\,v_{j}>}{\|v_{j}\|^{2}}\ast v_{j}\,,\;\;v_{i}>\smallskip\hspace{5cm}\]
\[\Rightarrow\;<v_{k}\,,\,v_{i}>\,\;\;\leq\;\;\,\sup< w_{k}\,\#\,- \mathop {\sum }\limits_{j=1 }\limits^{k-1}\frac{<w_{k}\,,\,v_{j}>}{\|v_{j}\|^{2}}\ast v_{j}\,,\;\;v_{i}>\smallskip\hspace{5cm}\]
\[=\,< w_{k}\,,\,v_{i}>\,+\,\sup<- \mathop {\sum }\limits_{j=1 }\limits^{k-1}\frac{<w_{k}\,,\,v_{j}>}{\|v_{j}\|^{2}}\ast v_{j}\,,\;\;v_{i}>\smallskip\hspace{5cm}\]
=$\,< w_{k}\,,\,v_{i}>\,-\,\frac{<w_{k}\,,\,v_{i}>}{\|v_{i}\|^{2}}<v_{i}\,,\,v_{i}>$\\
=$\,< w_{k}\,,\,v_{i}>\,-\,\frac{<w_{k}\,,\,v_{i}>}{\|v_{i}\|^{2}}\|v_{i}\|^{2}$\\
=$\,0.$\\
Therefore $<v_{k}\,,\,v_{i}>\;\;\,\leq 0\;\;\;\cdots\cdots\cdots (2)$\\
Again \[v_{k}\,\in\, w_{k}\,\#\,- \mathop {\sum }\limits_{j=1 }\limits^{k-1}\frac{<w_{k}\,,\,v_{j}>}{\|v_{j}\|^{2}}\ast v_{j}\smallskip\hspace{10cm}\]
 \[\Rightarrow(-1)\ast v_{k}\,\subseteq\, (-1)\ast( \,w_{k}\,\#\,- \mathop {\sum }\limits_{j=1 }\limits^{k-1}\frac{<w_{k}\,,\,v_{j}>}{\|v_{j}\|^{2}}\ast v_{j}\,)\smallskip\hspace{6cm}\]
\[\Rightarrow(-1)\ast v_{k}\,\subseteq\,(-1)\ast w_{k}\,\#\, \mathop {\sum }\limits_{j=1 }\limits^{k-1}\frac{<w_{k}\,,\,v_{j}>}{\|v_{j}\|^{2}}\ast v_{j}\smallskip\hspace{7cm}\]
\[=\,-1\ast w_{k}\,\#\, \mathop {\sum }\limits_{j=1 }\limits^{k-1}\frac{<w_{k}\,,\,v_{j}>}{\|v_{j}\|^{2}}\ast v_{j}\smallskip\hspace{3.5cm}\]
Therefore for $1\,\leq\,i\,\leq\,k-1$, we have\\
\[<-1\ast v_{k}\,,\,v_{i}>\,\;\;\subseteq\;\;\,<-1\ast w_{k}\,\#\, \mathop {\sum }\limits_{j=1 }\limits^{k-1}\frac{<w_{k}\,,\,v_{j}>}{\|v_{j}\|^{2}}\ast v_{j}\,,\;\;v_{i}>\smallskip\hspace{5cm}\]
\[\Rightarrow\,-< v_{k}\,,\,v_{i}>\,\;\;\leq\,\;\;\sup<-1\ast w_{k}\,\#\, \mathop {\sum }\limits_{j=1 }\limits^{k-1}\frac{<w_{k}\,,\,v_{j}>}{\|v_{j}\|^{2}}\ast v_{j}\,,\;\;v_{i}>\smallskip\hspace{5cm}\]
\[\smallskip\hspace{3.3cm}=\,\sup<-1\ast w_{k}\,,\,v_{i}>\,+\,\sup< \mathop {\sum }\limits_{j=1 }\limits^{k-1}\frac{<w_{k}\,,\,v_{j}>}{\|v_{j}\|^{2}}\ast v_{j}\,,\;\;v_{i}>\]\\
$\smallskip\hspace{3.3cm}=\,-<w_{k}\,,\,v_{i}>\,+\,\frac{<w_{k}\,,\,v_{i}>}{\|v_{i}\|^{2}}\|v_{i}\|^{2}$\\
$\smallskip\hspace{3.3cm}=\,0.$\\
Therefore $-<v_{k}\,,v_{i}>\,\leq\, 0\;\;\;\cdots\cdots\cdots(3)$\\
From (2) and (3), we get $<v_{k}\,,\,v_{i}>\,=\,0$ for $1\,\leq \,i\,\leq\,k-1$\\
Again by the induction hypothesis $S^{'}_{K-1}$ is orthogonal.\\
Hence $S^{'}_{K}$ is an orthogonal set of non-zero vectors.\\
We now show that $ span(S^{'}_{k})\,= \,span(S_{k})$\\
From $(1)$ we have \[v_{k}\,\in\, w_{k}\,\#\,- \mathop {\sum }\limits_{j=1 }\limits^{k-1}\frac{<w_{k}\,,\,v_{j}>}{\|v_{j}\|^{2}}\ast v_{j},\;\;\,2\,\leq\,k\,\leq\,n\]
\[\Rightarrow w_{k}\,\in\, v_{k}\,\#\, \mathop {\sum }\limits_{j=1 }\limits^{k-1}\frac{<w_{k}\,,\,v_{j}>}{\|v_{j}\|^{2}}\ast v_{j},\;\;\,2\,\leq\,k\,\leq\,n\]
Therefore $span(S_{k})\,\subseteq\,span(S^{'}_{k})\;\;\;\cdots\cdots\cdots (4)$\\
Next let $\alpha\,\in\,span(S^{'}_{k})$, then there exists $\lambda_{1},\,\lambda_{2},\,\cdots,\,\lambda_{k}\,\in\,\mathbf{R}$ such that\\
$ \alpha\,\in\,\lambda_{1}\ast v_{1}\,\#\,\lambda_{2}\ast v_{2}\,\#\,.\,.\,.\,\#\,\lambda_{k}\ast v_{k}$\\
At first by $(1)$ we replace $v_{k}$ with \[w_{k}\,\#\,- \mathop {\sum }\limits_{j=1 }\limits^{k-1}\frac{<w_{k}\,,\,v_{j}>}{\|v_{j}\|^{2}}\ast v_{j}\]
Then there exists $\delta_{1},\,\delta_{2},\,\cdots,\,\delta_{k}\,\in \mathbf{R}$ such that \\
$ \alpha\,\in\,\delta_{1}\ast v_{1}\,\#\,\delta_{2}\ast v_{2}\,\#\,.\,.\,.\,\#\,\delta_{k-1}\ast v_{k-1}\,\#\,\delta_{k}\ast w_{k}$.\\
In the same way, then we replace $v_{k-1}$ with \[w_{k-1}\,\#\,- \mathop {\sum }\limits_{j=1 }\limits^{k-2}\frac{<w_{k-1}\,,\,v_{j}>}{\|v_{j}\|^{2}}\ast v_{j}\;\;\;[ by \,(1)]\]
proceeding in above repeated ways, lastly we find the scalars \\
$\eta_{1},\,\eta_{2},\,\cdots,\,\eta_{k}\,\in\,\mathbf{R}$ such that
$\; \alpha\,\in\,\eta_{1}\ast w_{1}\,\#\,\eta_{2}\ast w_{2}\,\#\,.\,.\,.\,\#\,\eta_{k}\ast w_{k}$\\
 $\Rightarrow\alpha\,\in\,span(S_{k})$\\
 Therefore $span(S^{'}_{k})\,\subseteq\,span(S_{k})\;\;\;\cdots\cdots\cdots (5)$\\
 So from $(4)$ and $(5)$ we get $span(S^{'}_{k})\,=\,span(S_{k})$.
 Hence by induction $S^{'}$ is an orthogonal set of non-zero vectors such that \\$span(S^{'})\,=\,span(S)$\\
 This completes the proof.\\
 \smallskip\hspace{.5cm} The construction of the set $S^{'}$ by the use of theorem $\ref{t1}$ is called the Gram-Schmidt process.\\
 \begin{note}
 every finite dimensional innerproduct hyperspace possesses an orthogonal basis.
\end{note}
\textbf{Acknowledgements}\\
  The authors are grateful to the referees for their valuable suggestions in rewriting the paper in the present form.\\


\begin{thebibliography}{0}
\bibitem{Ameri}  Ameri, R.
\emph{Fuzzy Hypervector Spaces Over Valued Field},
\bibitem{Corsini} Corsini, p.  Leoreanu, V.
\emph{Applications of Hyperstructure Theory}, Kluwer Academic Publishers, Dordrecht, Hardbound, (2003).
\bibitem{Krasner} Krasner, M.
\emph{A class of hyperrings and hyperfields}, Intern. J.Math. and Math. Sci .,Vol6, no.2, (1983),307-312,.
\bibitem{Marty1}
\emph{Marty, F. Sur une g´en´eralisation de la notion de groupe.}, In 8`eme congr`es des Math´ematiciens
Scandinaves, Stockholm, pages 45-49, (1934).
\bibitem{Marty2}
\emph{Marty, F. Rˆole de la notion d'hypergroupe dans l'´etude des groupes non ab´eliens}, Comptes Renclus Acad. Sci. Paris Math, 201, 636-638,$(\,1935\,)$.
\bibitem{Marty3}
\emph{Marty, F. Sur les groupes et hypergroupes attach´es `a une fraction rationnelle}, Ann. Sci. de l' Ecole Norm. Sup., (3) 53, 82-123, $(\,1936\,)$.
\bibitem{Nakassis}Nakassis A.
\emph{Expository and survey article recent results in hyperring and hyperfield theory}, internet. J. Math. and Math. Sci.,11(2),(1988),209 - 220.
\bibitem{Roy}Roy S, Samanta T, K.
\emph{A note on Hypervector Spaces}, (accepted in Discussiones Mathematicae - General Algebra and Applications)
\bibitem{Tallini1}Scafati-Tallini, M.
\emph{ Hypervector spaces}. 4th AHA, World Scientific, Xanthi, Greece, (1991), 167-174.
\bibitem{Tallini2}Scafati-Tallini, M.
 \emph{Weak hypervector spaces and norms in such spaces}. Proceedings of the 5th International Congress on AHA and Appl.(1993), Hadronic Press., Jasi,Rumania, (1994), 199-206.
\bibitem{Tallini3}Scafati-Tallini, M.,
\emph{Characterization of remarkable hypervector spaces}., Proceedings of the 8th International Congress on AHA and Appl., Samothraki, Greece, (2003), 231-238.

\end{thebibliography}
\end{document}